\documentclass{article}

\usepackage[includeheadfoot,
            bindingoffset=0mm,
            inner   = 25mm,
            top     = 20mm,
            outer   = 25mm,
            bottom  = 10mm,
            paperwidth = 210mm,
            paperheight = 297mm,
            ]{geometry}
\usepackage[margin={2.0cm,0cm},oneside,labelfont={sf,bf},singlelinecheck=false]{caption}

\usepackage{graphicx}

\usepackage[fleqn]{amsmath}
\usepackage{accents}
\usepackage{amssymb}
\usepackage{amsmath}

\usepackage{enumitem}
\usepackage{cite}
\usepackage[perpage]{footmisc}
\usepackage{hyphenat}
\usepackage[english]{babel}
\usepackage[section]{placeins}
\usepackage{color}
\usepackage{upgreek}
\usepackage{listings}

\definecolor{codegreen}{rgb}{0,0.6,0}
\definecolor{codeblue}{rgb}{0,0,0.8}
\definecolor{codegrey}{rgb}{0.5,0.5,0.5}
\usepackage{amsthm}
\theoremstyle{definition}

\newcommand\undermat[2]{\makebox[0pt][l]{$\smash{\underbrace{\phantom{\begin{matrix}#2\end{matrix}}}_{#1}}$}#2}

\usepackage{titlesec}
\titleformat{\section}[hang]{\Large\bfseries\raggedright\sffamily}{\thesection}{1em}{}
\titleformat{\subsection}[hang]{\large\bfseries\raggedright\sffamily}{\thesubsection}{1em}{}
\titleformat{\subsubsection}[hang]{\normalsize\bfseries\raggedright\sffamily}{\thesubsubsection}{1em}{}
\usepackage{abstract}

\usepackage{authblk}

\newcommand{\transp}{\ensuremath{ ^\mathrm{T} }}

\begin{document}

\title{ \huge\bfseries\sffamily A deep machine learning algorithm for construction of the Kolmogorov-Arnold representation }

\author[1]{A. Polar}
\author[2]{M. Poluektov}
\affil[1]{Independent Software Consultant, Duluth, GA, USA}
\affil[2]{International Institute for Nanocomposites Manufacturing, WMG, University of Warwick, Coventry CV4 7AL, UK}

\date{ \huge\normalfont\sffamily DRAFT: \today }

\maketitle

\setlength{\absleftindent}{2.0cm}
\setlength{\absrightindent}{2.0cm}
\setlength{\absparindent}{0em}
\begin{abstract}
The Kolmogorov-Arnold representation is a proven adequate replacement of a continuous multivariate function by an hierarchical structure of multiple functions of one variable. The proven existence of such representation inspired many researchers to search for a practical way of its construction, since such model answers the needs of machine learning. This article shows that the Kolmogorov-Arnold representation is not only a composition of functions but also a particular case of a tree of the discrete Urysohn operators. The article introduces new, quick and computationally stable algorithm for constructing of such Urysohn trees. Besides continuous multivariate functions, the suggested algorithm covers the cases with quantised inputs and combination of quantised and continuous inputs. The article also contains multiple results of testing of the suggested algorithm on publicly available datasets, used also by other researchers for benchmarking. discrete Urysohn operator, classification trees.
\end{abstract}

\section{Introduction}
\label{sec:intro}

The Kolmogorov-Arnold representation \cite{Kolmogorov1956} of a continuous multivariate function is a decomposition of the function into a structure of inner and outer functions of a single variable. More precisely, function $F: \mathbb{R}^m \to \mathbb{R} \in C\left(\left[0,1\right]^m\right)$ can be represented as
\begin{equation}
  F\left( x_1, x_2, \ldots, x_m \right) = \sum_{k=1}^{2m+1} \varPhi^k \left( \sum_{j=1}^{m} f^{kj} \left(x_j\right) \right) ,
  \label{eq:Kolmogorov} 
\end{equation}
where $f^{kj}: \left[0,1\right] \to \left[0,1\right] \in C\left[0,1\right]$ and $\varPhi^k: \mathbb{R} \to \mathbb{R} \in C\left(\mathbb{R}\right)$. The original article \cite{Kolmogorov1956} only states the existence of such representation and does not provide a method for its construction.

The further research of this representation can be conventionally divided into generic ways of model reduction \cite{Sprecher1965,Lorentz1966,Sprecher1972,Lorentz1996} and practical ways of construction of the involved functions after picking one of the existing reduced forms 
\cite{Igelnik2003,Coppejans2004,Wasserman2006,Actor2017,Actor2018,Montanelli2020}. Furthermore, the relation of representation \eqref{eq:Kolmogorov} to neural networks has been noticed and investigated by several researchers \cite{HechtNielsen1987,Lippmann1987,Kurkova1992,Sprecher1993}.

This article approaches the construction (or \emph{identification}) of representation \eqref{eq:Kolmogorov} from a practical perspective. While a generic input-output relationship can be described by a continuous multivariate function, it is possible that for a particular dataset, a reduced form of representation \eqref{eq:Kolmogorov} may be sufficient. Furthermore, it may also be sufficient to represent the underling functions of representation \eqref{eq:Kolmogorov} up to a certain numerical accuracy, as long as it can be controlled by the number of unknowns. The aim of this article is to propose an efficient identification algorithm, which incorporates such control over the identified model.

A previously unnoticed (to the best knowledge of the authors) property of representation \eqref{eq:Kolmogorov} is the connection to the discrete Urysohn operator \cite{Krylov1979}, which transforms a sequence of scalars into another scaler. More precisely, $U: \mathbb{R}^m \to \mathbb{R}$ is given by
\begin{equation}
  U\left( x_1, x_2, \ldots, x_m \right) = \sum_{j=1}^{m} g^j \left(x_j\right) ,
  \label{eq:UrysDiscr} 
\end{equation} 
where functions $g^j: \mathbb{R} \to \mathbb{R}$ are functions of one variable. The only formal distinction between the functions of the discrete Urysohn operator and the Kolmogorov-Arnold representation is that $g^j$ are subject to fewer restrictions, e.g. they might have a discontinuity of the first kind. After introduction of auxiliary intermediate parameters $\phi_k$, it becomes obvious that the Kolmogorov-Arnold representation is also a \emph{tree} of the discrete Urysohn operators with a single root operator and $\left(2m + 1\right)$ branch operators:
\begin{equation}
  F\left( x_1, x_2, \ldots, x_m \right) = \sum_{k=1}^{2m+1} \varPhi^k \left( \phi_k \right) , \quad\quad
  \phi_k = \sum_{j=1}^{m} f^{kj} \left(x_j\right) .
  \label{eq:UrysDiscrRoot} 
\end{equation}

The identification algorithm, suggested in this article, for the Kolmogorov-Arnold representation (or a tree of the discrete Urysohn operators) is based on recently published research \cite{Poluektov2019} by the same authors on non-parametric identification of an individual discrete Urysohn operator for given input-output data. The novelty of this publication is in advancement of previously-published method \cite{Poluektov2019} to the case when the discrete Urysohn operators are arranged in a chain with unobserved intermediate values, as in equation \eqref{eq:UrysDiscrRoot}. The suggested method uses each Urysohn operator as a single element in the identification and updates all functions of each operator at one step (i.e. the functions are not treated individually), which significantly expedites the entire identification process. The article targets practical aspects rather than theory and is backed up by downloadable and reproducible tests.

Application of representation \eqref{eq:Kolmogorov} to modelling of physical or social systems introduces new aspects, not considered in the original work \cite{Kolmogorov1956}. For a physical system, the data may contain measurement noise, while for a social system, the data may even be stochastic. For example, two individuals with identical demographic parameters (the inputs) can make different decisions regarding a purchase of goods or services (the binary output). When representation \eqref{eq:Kolmogorov} is used for modelling of approximate data, it is interpreted as a model of a particular multivariate function, which minimises the error between calculated output $F$ and the real (measured) output for the unseen data.  

\section{Identification of the single Urysohn operator}
\label{sec:basic_method}

The identification of a tree of Urysohn operators is based on the identification of a single operator for input-output data. This method is published in \cite{Poluektov2019}, where a detailed research of both quantised and continuous discrete Urysohn operators, including ones with multiple vector inputs, is given. It must be mentioned, that \cite{Poluektov2019} discusses the Urysohn operators from a different perspective --- modelling dynamic systems; therefore, the notation there is slightly different to this article. For convenience of the readers, the basic concept is repeated in this section, which represents only a short digest, sufficient for understanding of the subsequent sections. 

\subsection{From linear regression to the Urysohn operator}
\label{sec:basic_method_explain}

At the first explanation step, the linear regression model is considered:
\begin{equation}
  \hat{z}_i = \sum_{j=1}^m w_j x_{j,i} ,
  \label{eq:Linear} 
\end{equation}
where $x_{j,i} \in \left[ x_{j,\mathrm{min}}, x_{j,\mathrm{max}} \right]$ is the $j$-th input of the $i$-th record, $\hat{z}_i$ is the calculated model output of the $i$-th record. Model \eqref{eq:Linear} is expanded into another linear model which has two points per given input (left $L_j$ and right $R_j$) instead of one weight coefficient $w_j$:
\begin{equation}
  \hat{z}_i = \sum_{j=1}^m \left( L_j \left(1 - p_{j,i}\right) + R_j p_{j,i} \right) , \quad\quad 
  p_{j,i} = \frac{x_{j,i} - x_{j,\mathrm{min}}}{x_{j,\mathrm{max}} - x_{j,\mathrm{min}}} .
  \label{eq:Linear2} 
\end{equation}
A large set of input-output records for model \eqref{eq:Linear2} forms a system of linear algebraic equations 
with unknown vector-column
\begin{equation}
  V = \begin{bmatrix}
    L_1 & L_2 & \ldots & L_m & R_1 & R_2 & \ldots & R_m
  \end{bmatrix}\transp
  \label{Vector}
\end{equation}
and a matrix with rows $P_i$ built out of input values
\begin{equation}
  P_i = \begin{bmatrix}
    1-p_{1,i} & 1-p_{2,i} & \ldots & 1-p_{m,i} & p_{1,i} & p_{2,i} & \ldots & p_{m,i} 
  \end{bmatrix} .
  \label{Vector2}
\end{equation}
The matrix is always singular independently of data, but the solution, can be obtained as an approximation by minimising the norm of $V$, using model equations \eqref{eq:Linear2} as constraints. There are different ways of finding the constrained minimum; the authors have chosen the projection descent method \cite{Kaczmarz1937,Tewarson1969,Faddeev1981,Hykin1996}. In this method, $V$ is iteratively changed for each input-output record independently. It is convenient to denote $V$ at iteration $i$ as $V_i$.

In the projection descent method, the initial approximation $V_1$ is assigned, which is the all-zero vector-column for minimising $\left|V\right|$. At each iteration, for each new data record, the following modification is performed:
\begin{equation}
  V_{i+1} = V_i + \alpha \frac{z_i - P_i V_i}{\left|P_i\right|^2} {P_i}\transp ,
  \label{projection}
\end{equation}
where $z_i$ is the real (recorded) output of the $i$-th record and parameter $\alpha \in \left(0,2\right)$ provides error filtering and controls the convergence rate. Here, $\left|P_i\right|^2 \in \left[m/2,m\right]$ due to the definition of $p_{j,i}$. Numerator $D_i = z_i - P_i V_i$ in equation \eqref{projection} is the residual. The matter of the projection descent method \cite{Kaczmarz1937} is the navigation of point $V_i$ to the solution by projecting it from one hyperplane to another. The magnitude of $\alpha D_i$ controls the distance and the sign of $D_i$ controls the direction. Having the correct sign is crucial, as the wrong sign navigates the projected point away from the solution. 

Now, the model based on the discrete Urysohn operator \eqref{eq:UrysDiscr} is considered and, for convenience, it is rewritten as
\begin{equation}
  \hat{z}_i = \sum_{j=1}^{m} g^j \left(x_{j,i}\right) ,
  \label{eq:UryModel}
\end{equation} 
where $x_{j,i}$ and $\hat{z}_i$ are the inputs and the calculated model output of the $i$-th record, respectively. The next key step is that functions $g^j$ are sought in the class of piecewise-linear functions. The nodal values for each function $g^j$, i.e. the function values where it changes the slope, are denoted as $G^j_k$. In this case, the projection descent method is applicable with a few changes. Estimated vector-column $V$ now has a block structure, where all nodal values $G^j_k$ are written in the sequential order:
\begin{equation*}
  V = \begin{bmatrix}
    G^1_1 & G^1_2 & \ldots & G^1_{n_1} & G^2_1 & G^2_2 & \ldots & G^2_{n_2} & \ldots & G^m_1 & G^m_2 & \ldots & G^m_{n_m}
  \end{bmatrix}\transp ,
\end{equation*}
where $n_j$ is the number of the nodes of function $g^j$. Argument $x_{j,i}$ will always\footnote{If $x_{j,i}$ falls exactly onto a nodal position, any of the adjacent linear segments can be taken. In either case, one of $\left(1-\psi_{j,i}\right)$ and $\psi_{j,i}$ will be $0$, while another will be $1$.} fall into only one linear segment for each function $g^j$, hence, relative distance $\psi_{j,i}$ between the beginning of the linear segment and $x_{j,i}$ can be introduced. Vector-row $P_i$ also has a corresponding block structure with only one pair of non-zero elements per block, which are $\left(1-\psi_{j,i}\right)$ and $\psi_{j,i}$ for the argument in this block,
\begin{equation*}
  P_i = \begin{bmatrix}
    P^1_i & P^2_i & \ldots & P^m_i 
  \end{bmatrix} , \quad\quad
  P^j_i = \begin{bmatrix}
    \undermat{q_{j,i}-1}{0 & \ldots & 0} & 1-\psi_{j,i} & \psi_{j,i} & \undermat{n_j-q_{j,i}-1}{0 & \ldots & 0} 
  \end{bmatrix}
  \vphantom{\underbrace{\begin{matrix}0 & \ldots & 0\end{matrix}}_{q_{j,i}}} ,
\end{equation*} 
where $q_{j,i}$ is the number of the linear segment, into which $x_{j,i}$ falls. The matrix of the system becomes sparse, compared to the considered above case, and the projection descend method converges faster, since the adjacent rows are either orthogonal or near-orthogonal.

\begin{figure}
  \begin{center}
    \includegraphics{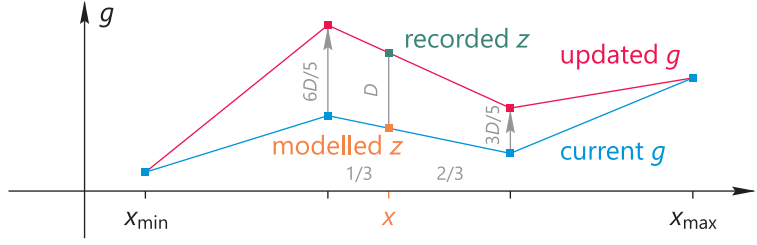}
  \end{center}
  \caption{A schematic illustration of one step of the iterative algorithm for identification of the discrete Urysohn operator. For clarity of the representation, the operator consisting of only one function is taken. The function is split into three linear segments. At the illustrated iteration, the input falls into the second segment and the nodes of that segment are modified according to the relative distances between the input and the nodes.}
  \label{fig:USurf2}
\end{figure}

In the computational implementation, matrices and vectors are not built --- the linear algebra notations and operations were used only for the explanation. The nodes that must be modified are identified by inputs $x_{j,i}$, the difference between the calculated model output and the actual output is determined and all involved nodal values are modified, as schematically illustrated in figure \ref{fig:USurf2}. The formal algorithm is given below.

\subsection{Formal algorithm for the Urysohn operator identification}
\label{sec:basic_method_algo}

First, functions $g^j$ must be written in a piecewise-linear form. For each function $g^j$, input interval $\left[ x_{j,\mathrm{min}}, x_{j,\mathrm{max}} \right]$ is divided into $\left(n_j-1\right)$ equal intervals; the rescaled inputs and their rounding to the nearest integer values are introduced as
\begin{equation}
  b_{j,i} = 1 + \left( n_j-1 \right) \frac{x_{j,i} - x_{j,\mathrm{min}}}{x_{j,\mathrm{max}} - x_{j,\mathrm{min}}} , \quad\quad
  q_{j,i} = \left\lfloor b_{j,i} \right\rfloor , \quad\quad
  r_{j,i} = \left\lceil b_{j,i} \right\rceil ,
\end{equation}
where $\lfloor \cdot \rfloor$ and $\lceil \cdot \rceil$ are the floor and the ceiling functions, respectively. Now, functions $g^j$ can be written as
\begin{equation}
  g^j \left( x_{j,i} \right) = \left( 1 - \psi_{j,i} \right) G^j\left[ q_{j,i} \right] +
  \psi_{j,i} G^j\left[ r_{j,i} \right] , \quad\quad
  \psi_{j,i} = b_{j,i} - q_{j,i} ,
\end{equation}
where $G^j$ is the vector-column containing the nodal values of function $g^j$, with indices shown in $\left[\cdot\right]$. Finally the following norm is introduced:
\begin{equation}
  \chi_i = \sum_{j=1}^{m} \left( \left( 1 - \psi_{j,i} \right)^2 + {\psi_{j,i}}^2 \right) .
\end{equation}
The algorithm starts with all $G^j$ being the all-zero vector-columns. The algorithm then proceeds iteratively, with one iteration consisting of the following steps: 
\begin{enumerate}
  \item Calculate model output $\hat{z}_i$ based on actual inputs $x_{j,i}$ and the current approximation of model parameters $G^j$;
  \item Calculate difference $D_i = z_i - \hat{z}_i$, where $z_i$ is the real (recorded) output and $\hat{z}_i$ is the calculated model output;
  \item Modify parameters $G^j$, such that $\alpha D_i \left( 1 - \psi_{j,i} \right) / \chi_i$ is added to parameters $G^j\left[q_{j,i}\right]$ and $\alpha D_i \psi_{j,i} / \chi_i$ is added to parameters $G^j\left[r_{j,i}\right]$ for each $j$.
\end{enumerate}
Here, parameter $\alpha \in \left(0, 2\right)$ is the noise-reduction parameter. The algorithm is considered to be converged, when $D_i$ becomes sufficiently small for sufficiently large number of iterations consecutively.

It must be noted that it is easy to write a generalisation of the above algorithm for the case when the distances between the neighbouring nodes of function $g^j$ may vary. Furthermore, the model and the algorithm, without any modifications, are also applicable for the quantised inputs. This case can be treated as if all inputs $x_{j,i}$ always fall onto the edges of the linear segments and only one nodal value per function is modified in the algorithm.

\section{Identification of an Urysohn tree}
\label{sec:Kolm}

As discussed in the introduction, the identification problem for the Urysohn tree consists in finding the unknown functions of the operators for a known set of input-output records. Furthermore, these functions are assumed to be piecewise linear, which means that the identification problem consists in finding the nodal values of the functions (i.e. the function values where the slope changes).

It is convenient to rewrite equation \eqref{eq:UrysDiscrRoot} and introduce the notation for the $i$-th input-output record in the discrete form:
\begin{equation}
  \hat{z}_i = \hat{z}_i\left( \phi_{k,i} \right) = \sum_{k=1}^{2m+1} \varPhi^k \left( \phi_{k,i} \right) , \quad\quad
  \phi_{k,i} = \sum_{j=1}^{m} f^{kj} \left(x_{j,i}\right) ,
  \label{eq:VKRoot} 
\end{equation}
where $x_{j,i}$ is the $j$-th input of the $i$-th record, $\hat{z}_i$ is the calculated model output of the $i$-th record. Obviously, if values $\phi_{k,i}$ are known, the identification problem for model \eqref{eq:VKRoot} reduces to the identification of multiple individual discrete Urysohn operators. Unfortunately, values $\phi_{k,i}$ are not only unobserved, but also exist only as auxiliary mathematical variables. 

Before presenting the identification steps, additional quantities must be introduced. Assuming that the current approximations for functions $\varPhi^k$ and $f^{kj}$ are given and the $i$-th data record is considered, the residual and the model error can be introduced as
\begin{equation*}
  R_i \left( \phi_{k,i} \right) = z_i - \hat{z}_i\left( \phi_{k,i} \right) , \quad\quad  
  E_i \left( \phi_{k,i} \right) = \left| R_i \left( \phi_{k,i} \right) \right| ,
\end{equation*}
where $z_i$ is the real (recorded) output of the $i$-th record. Next, the increments for auxiliary variables are defined as
\begin{equation}
  \Delta \phi_{k,i} = \frac{\mu R_i \left( \phi_{k,i} \right)}{\left(2m+1\right) T\left( {\varPhi^k}'\left( \phi_{k,i} \right) \right) } , \quad\quad
  T\left( \zeta \right) = \begin{cases} \zeta &\text{if } \left|\zeta\right| \geq \delta \\
    \delta & \text{if } 0 \leq \zeta < \delta \\
    -\delta & \text{if } -\delta < \zeta < 0 \end{cases} ,
  \label{eq:delta} 
\end{equation}
where $\mu \in \left( 0, 1 \right]$ is the noise-reduction parameter, function $T\left( \zeta \right)$ is introduced to account for the case of zero derivative of $\varPhi^k$ and $\delta$ is a small threshold value. It can be seen that if $\mu = 1$ and if functions $\varPhi^k$ are linear, with the derivatives greater than $\delta$, then the residual of updated auxiliary variables $R_i \left( \phi_{k,i} + \Delta \phi_{k,i} \right)$ becomes zero. In general, for piecewise-linear functions $\varPhi^k$ with the non-zero derivatives and small enough $\mu$, updating the auxiliary variables should decrease the model error\footnote{Except in cases when $\left( \phi_{k,i} + \Delta \phi_{k,i} \right)$ fall onto the nodal positions of $\varPhi^k$, where the sign of its derivative changes from negative to positive.}:
\begin{equation}
  E_i \left( \phi_{k,i} + \Delta \phi_{k,i} \right) < E_i \left( \phi_{k,i} \right) .
  \label{eq:ERR2} 
\end{equation}
However, for a prescribed $\mu$, inequality \eqref{eq:ERR2} might not hold.

\subsection{Record-by-record descent}

The basic idea behind the identification algorithm presented here is to try to change auxiliary variables $\phi_{k,i}$ by $\Delta \phi_{k,i}$, verify whether it decreases the modelling error and, if it does, then assume $\phi_{k,i} + \Delta \phi_{k,i}$ to be the new (known) auxiliary variables, based on which all operators are updated. It is useful to remind the utilised terminology --- the branch Urysohn operators contain functions $f^{kj}$ and the root Urysohn operator contains functions $\varPhi^k$.

The formal steps of algorithm are as follows:
\begin{enumerate}
  \item Make initial approximations of the Urysohn operators (see section \ref{sec:init}); 
  \item Take one input-output record, calculate auxiliary variables $\phi_{k,i}$, model output $\hat{z}_i$, model error $E_i\left( \phi_{k,i} \right)$ and increments $\Delta \phi_{k,i}$, given the current approximations of the Urysohn operators; \label{st:calc}
  \item Calculate new auxiliary variables $\phi_{k,i} + \Delta \phi_{k,i}$ and calculate model error $E_i \left( \phi_{k,i} + \Delta \phi_{k,i} \right)$ corresponding to the new auxiliary variables; \label{st:new}
  \item If inequality \eqref{eq:ERR2} does not hold, then go to step \ref{st:calc}, otherwise proceed;
  \item Update all branch Urysohn operators based on the new auxiliary variables obtained at step \ref{st:new} (see section \ref{sec:basic_method_algo} for the single operator);
  \item Change the nodal positions of the root Urysohn operator, if the new auxiliary variables are outside of its domain (see below). \label{st:updNode}
  \item Update the root Urysohn operator based on the new auxiliary variables obtained at step \ref{st:new}; \label{st:root}
  \item Check for convergence, if not converged, then go to step \ref{st:calc}. 
\end{enumerate}

During the update of the auxiliary variables at step \ref{st:new}, it can happen that some of new auxiliary variables $\phi_{k,i} + \Delta \phi_{k,i}$ are outside of the domain of the functions of the root Urysohn operator. In this case, the nodal positions of corresponding functions $\varPhi^k$ must be updated (step \ref{st:updNode}) --- new domain limits are calculated and the nodal positions are redistributed to cover the new domain. The function values at the new nodal positions are recalculated using either interpolation (if the new node is inside the old domain) or extrapolation (if the new node is outside the old domain).

The algorithm is specifically presented for the Kolmogorov-Arnold representation, as such representation is sufficient to describe any continuous multivariate function; however, it is easy to construct the generalisation of the algorithm for the case when the tree has arbitrary many layers of the discrete Urysohn operators. In this case, at first, the increments for all auxiliary variables are found from top to bottom (i.e. ones that are closer to the output --- first), then all Urysohn operators are updated from bottom to top.

All conducted numerical experiments (see section \ref{sec:simulation}) confirmed quick convergence; however, at this moment, it is not backed up by a strict theoretical proof.

\subsection{Initial approximation} 
\label{sec:init} 

The identification of a single discrete Urysohn operator converges independently of an initial approximation. When all initial values are zeros, it converges to the solution with the minimum norm. The Urysohn tree, however, needs a specific initialisation. The simplest, but still an effective way is to assign random auxiliary variables $\phi_{k,i}$ and then update all branch operators using actual inputs $x_{j,i}$. Afterwards, calculate new auxiliary variables $\phi_{k,i}$ and, based on them, update the root operator for real (recorded) outputs $z_i$.

In the case if all branch operators are taken to be identical, the identification is much less efficient; thus it must be avoided. The reason for this is that at each identification step, all branch operators are updated in the same way and, if they all are initially identical, they stay identical, which is the same as having only one addend in the model (i.e. only one $\varPhi^k$).

\section{Numerical simulations and tests}
\label{sec:simulation} 

\subsection{Goals, accuracy metrics and overfitting}

The numerical simulations of this section pursue multiple goals, primary of which is comparing the proposed approach to the machine learning results published by other software companies. It must be emphasised that this is not a competition in accuracy, because data modelling is a special type of art, where a model can often be fine-tuned for a considerable time, until a previously published benchmark is beaten. Unfortunately, the fact that a particular method is more accurate for particular dataset than another method cannot be used to make any meaningful conclusions. The goal of the comparison presented here is to make sure that the accuracy is somehow near the level reported by other researchers or companies dealing with data modelling. 

The accuracy metrics are computed for the so-called unseen data or data that has not been used in the training process. The reported metrics are the Pearson correlation coefficient between actual $z_i$ and modelled $\hat{z}_i$ outputs and the normalised root mean square error (RMSE), defined as
\begin{equation}
  \bar{E}_\mathrm{RMSE} = \frac{1}{z_\mathrm{max}-z_\mathrm{min}}\sqrt{ \frac{1}{N}\sum_{j=1}^N \left(z_j - \hat{z_j}\right)^2 } ,
  \label{RMSE}
\end{equation}
where $z_\mathrm{min}$ and $z_\mathrm{max}$ are the minimum and the maximum values of the output, respectively. Furthermore, the error bounds are reported for each quantity as the $95\%$ confidence interval.

The overfitting for the Kolmogorov-Arnold representation can be easily prevented by a model reduction. According to the model for $m$ input parameters, the representation must have $2m + 1$ terms with $m$ functions each. As it has been found by testing, such large number of functions may not be required. Furthermore, the number of the linear blocks of each involved function can be varied and, in some cases, can be even as low as one (i.e. two points per function). Similarly to other commonly-used types of models, the reduced models give higher errors on the training dataset and lower errors on the validation dataset. The trade-off in accuracy and the model reductions are also demonstrated and explained in the performed tests. 

\subsection{Multivariate function}

The input-output records of physically existing systems may not be used for the assessment of the efficiency of the modelling concept because data may have hidden faults, such as omitted or ignored inputs, significant statistical differences between the training and the validation data, large experimental noise, etc. Therefore, the first test is conducted for generated data --- the ideal scenario, for which the method should give a $100\%$-accurate model.

The simulated system has $5$ inputs and $1$ output. It is computed by the following formula:
\begin{equation}
  z = \frac{\left|{\sin\left(x_2\right)}^{x_1}-\exp\left(-x_3\right)\right|}{x_4} + x_5 \cos\left(x_5\right) ,
  \label{formula}
\end{equation}
which has no physical meaning and is chosen only as a challenging expression. The inputs are generated randomly according to the uniform distribution within the following ranges: $x_1\in\left[0,0.99\right]$, $x_2\in\left[0,1.55\right]$, $x_3\in\left[1,1.49\right]$, $x_4\in\left[0.4,1.39\right]$, $x_5\in\left[0,0.49\right]$. The output range for equation \eqref{formula} with the given input limits is $z\in\left[0,2.37\right]$. The number of records has been taken to be $4000$.

Prior to applying the Kolmogorov-Arnold representation, the basic linear regression model, equation \eqref{eq:Linear}, has been tested. Obviously, equation \eqref{formula} is far from being linear and the linear regression has failed, as expected. The entire dataset has been used for the training and the Pearson correlation coefficient of $P = 0.88$ has been observed. 

The next test has been the single Urysohn operator, equation \eqref{eq:UryModel}. Again, the entire dataset has been used for the training. The Pearson correlation coefficient of $P = 0.93$ has been observed. The single Urysohn model introduces nonlinearities, therefore, has a higher accuracy than the linear model. However, it cannot reach the ideal accuracy, since each input value in equation \eqref{eq:UryModel}, after being transformed by a function, makes an additive contribution to the output, while equation \eqref{formula} contains inputs also in products.  

The Kolmogorov-Arnold model has been tested only on the unseen data, using the $10$-fold cross-validation method. After the records have been generated, $90\%$ of all records have been used for the training and the remaining $10\%$ for the validation. The test has been repeated $10$ times, each time with a different validation segment, such that each record has been used $9$ times for the training and $1$ time for the validation. For constructing the model, $11$ addends in the Kolmogorov-Arnold model (i.e. $2m + 1$, where $m$ is the number of the inputs) and $10$ nodal points per function have been used.

Based on $5$ consecutive executions, the normalised RMSE of $\bar{E}_\mathrm{RMSE} = 0.0203 \pm 0.0023$ and the Pearson correlation coefficient of $P = 0.9935 \pm 0.0014$ have been observed. This is the expected performance of the Kolmogorov-Arnold model, approaching the ideal accuracy. The discrepancy can be attributed to a finite number of nodal points per function and a finite number of records used for the identification. Therefore, it can be concluded that the proposed method passes this particular test of modelling the input-output data generated by the non-linear multivariate function.

To study the effects of the model reduction, the number of addends in the Kolmogorov-Arnold model has been varied. For a single addend, the Pearson correlation coefficient of $P = 0.88$ has been observed, for two addends --- $P = 0.93$, for three addends --- $P = 0.96$. Finally, for $6$ addends, $P$ becomes comparable to that of the full model, which means that for this particular example, it is not necessary to use all $11$ addedns in the Kolmogorov-Arnold model for achieving a high accuracy. This is an illustration of how for many real life datasets, the model can be reduced and the researchers can vary the complexity by moving from the most simple to the full one. The source code and the data are available online\footnote{http://ezcodesample.com/naf/reallife5.html}. The code is reusable --- the generated data can be replaced by records for a physical system. 

\subsection{Airfoil self-noise}

The next example is a comparison of the proposed approach to the data science and machine learning platform Neural Designer\footnote{https://www.neuraldesigner.com/} developed by Artelnics. The chosen dataset is the ``Airfoil self-noise'', which can be found in UCI Machine Learning Repository \cite{Dua:2019}. The dataset is experimentally-obtained from a series of aerodynamic and acoustic tests of two- and three-dimensional airfoil blade sections conducted in an anechoic wind tunnel. The details of the Neural Designer modelling techniques can be found in the online help\footnote{https://www.neuraldesigner.com/learning/examples/airfoil-self-noise-prediction}. In this test, the dataset, consisting of $1503$ records, is randomly split into $60\%$ training subset, $20\%$ selection subset and $20\%$ testing subset. The model is trained multiple times on the training subset, tested on the unseen selection subset and the result with the best metric is then applied to the testing subset. The company reports that they have achieved the Pearson correlation coefficient of $P = 0.952$.

The modelled system has $5$ inputs and $1$ output. For constructing the model, $11$ addends in the Kolmogorov-Arnold model (i.e. exactly according to equation \eqref{eq:Kolmogorov}) and $15$ nodal points per function have been used. To illustrate the dataset, the first three records are shown below:
\begin{align*}
  &800;0;0.3048;71.3;0.00266337;126.201 \\
  &1000;0;0.3048;71.3;0.00266337;125.201 \\
  &1250;0;0.3048;71.3;0.00266337;125.951 .
\end{align*}
The authors followed the same training-testing path as reported by Artelnics --- the programme repeated training $10$ times, chose the best model according to the metric obtained on the selection subset and applied this model to the testing data. Based on $5$ consecutive executions, the Pearson correlation coefficient of $P = 0.9506 \pm 0.0049$ has been observed. This means that for the considered experimental dataset, the accuaracy of the proposed modelling and identification technique is at the same level as the accuary of the machine learning software Neural Designer. The source code and the data are available online\footnote{http://ezcodesample.com/naf/reallife3.html}.

\subsection{Mushroom classification}

The third example is a classification problem. This dataset has been publicly available for a long time and there are many software companies and individuals who tried modelling it. The input consists of $22$ mushroom properties and the output is a quality of the mushroom which is ether ``edible'' or ``poisonous''. The dataset ``Mushroom'' consists of $8124$ records and can be found in UCI Machine Learning Repository \cite{Dua:2019}. To illustrate the dataset, the first three records are shown below:
\begin{align*}
&\mathrm{p,x,s,n,t,p,f,c,n,k,e,e,s,s,w,w,p,w,o,p,k,s,u} \\
&\mathrm{e,x,s,y,t,a,f,c,b,k,e,c,s,s,w,w,p,w,o,p,n,n,g} \\
&\mathrm{e,b,s,w,t,l,f,c,b,n,e,c,s,s,w,w,p,w,o,p,n,n,m} \, .
\end{align*}
The first column contains the output --- ``edible'' or ``poisonous'', while other columns contain the mushroom properties.

Since the inputs are represented symbolically, some data preprocessing is required. Although the Kolmogorov-Arnold model relies on continuous functions, it can work perfectly with the quantised input data, as discussed in section \ref{sec:basic_method_algo}. For the quantised data, e.g. sequential integers, the nodes of the functions are arranged such that a particular input always coincide with a node. The functions then are assumed to be linear between the nodes, as for the non-quantised case. In the implementation, the symbolic inputs were transformed into the sequential integers, for example, the second column of the dataset contains elements of either of the following types: `b', `c', `f', `k', `s', `x', which were replaced by integers from $1$ to $6$. The symbolic outputs `e' and `p' were replaced by $1$ and $-1$, respectively.

For a comparison to a third party, the Microsoft ML.NET library\footnote{https://dotnet.microsoft.com/apps/machinelearning-ai/ml-dotnet} has been used. The accuracy for the library using the $10$-fold cross-validation method for the classification tree model was $57$ errors for entire dataset of $8124$ records. In a single step of the $10$-fold cross-validation method, $90\%$ of all records are used for the training and the remaining $10\%$ for the validation. Such step is repeated $10$ times, each time with a different validation segment. The execution time for the Microsoft ML.NET library was $63$ seconds.

The same $10$-fold cross-validation approach has been used for a single quantised Urysohn operator. For $10$ consecutive executions of the programme, the average number of mistakenly identified records was $3.60\pm0.37$. The execution time of the code was $2$ seconds. This result shows that even a single quantised Urysohn operator has outstanding descriptive capabilitites and can outperform established techniques, both in accuracy and execution time. The source code and the data are available online\footnote{http://ezcodesample.com/naf/reallife4.html}.

\subsection{Electronic non-linear dynamic system}

This test has been conducted for input-output recordings of an electronic non-linear dynamic system, which is assembled on a circuit board. The system is of a Wiener-Hammerstein type and contains a static nonlinearity that is sandwiched between two linear time-invariant blocks. Dataset ``Wiener-Hammerstein System (2009)'' \cite{Schoukens2009} has been taken from the website\footnote{http://www.nonlinearbenchmark.org/} specifically aimed at benchmarking of non-linear dynamic models, where the detailed description of the tested system and the data can be found. Unfortunately, the website deliberately does not provide the accuracy that they have achieved themselves, motivating this by saying that the ``benchmark is not intended as a competition'' \cite{Schoukens2009}.

Since the considered object is a single-input single-output (SISO) dynamic system, its inputs and outputs are arrays of numbers of size $1$ by $N$. However, an output $z_i$ of a dynamic system depends not only on $x_i$, but also on a finite fragment of preceding inputs, thus $z_i = z_i\left(x_i, x_{i-1}, x_{i-2}, \ldots, x_{i-m+1}\right)$. Therefore, for the identification of the Kolmogorov-Arnold model, the data has been rearranged accordingly, such that each record represents $m$ inputs and $1$ output.

The entire dataset of $N=188000$ records has been split equally into the training and the validation subsets, according to recommendation of the website. For the considered example, the length of the input fragment that defines the output has been determined to be $m=35$ elements. In this case, the full Kolmogorov-Arnold representation should have $71$ addends. However, only $4$ addends appeared to be sufficient and the number of the linear blocks in each function has been taken to be $15$. In addition to this, the linear regression model, equation \eqref{eq:Linear}, and the single Urysohn model, equation \eqref{eq:UryModel}, have been tested as well.

For the Kolmogorov-Arnold representation, the normalised RMSE of $\bar{E}_\mathrm{RMSE} = 0.0150 \pm 0.0016$ has been observed, while for the single Urysohn model $\bar{E}_\mathrm{RMSE} = 0.0152\%$ and for the linear model $\bar{E}_\mathrm{RMSE} = 0.032\%$ have been calculated. The difference between the Kolmogorov-Arnold model and the single Urysohn model is insignificant, since when the error is as low as $1.5\%$, there is almost no room for improvement, especially when the unseen data is of a physical system, with errors in the measurement. Also, the single Urysohn model is a particular case of the Kolmogorov-Arnold model (i.e. one branch and a linear outer function) and it is recommended as a starting point for modelling. Overall, such errors are usually considered to be accepatable for modelling dynamic systems. The source code and the data are available online\footnote{http://ezcodesample.com/naf/reallife1.html}.

\subsection{Bank churn}

The finial test is a model of the social system, where individuals make decisions regarding a subscription to a bank service. The data are individual features and the output is a human decision. This dataset is taken from the Neural Designer website\footnote{https://www.neuraldesigner.com/learning/examples/bank-churn}. The file contains the decisions regarding the bank clients and some of their personal data (provided anonymously), used as inputs. To illustrate the dataset, the first three records are shown below:
\begin{align*}
  &15634602;619;\mathrm{France;Female};42;2;0;1;1;1;101348.88;1 \\
  &15647311;608;\mathrm{Spain;Female};41;1;83807.86;1;0;1;112542.58;0 \\
  &15619304;502;\mathrm{France;Female};42;8;159660.8;3;1;0;113931.57;1 .
\end{align*}
The first number is the ID, which is ignored for modelling; the last binary value is YES/NO model output; the other parameters are either quantised, such as the country and the gender, or continuous, such as the balance and the salary. The total number of records is $10000$. 

The authors of Neural Designer split the dataset into $60\%$ training subset, $20\%$ selection subset and $20\%$ testing subset. The training is conducted multiple times, the best model is chosen based on the selection subset and is applied to the testing subset. They report the percentage of the correctly modelled records of $78.9\%$ on the testing subset.

For the modelling approach poposed in this article, the branch Urysohn operators have been built as having both continuous and quantised inner functions. The number of addends in the Kolmogorov-Arnold model have been taken to be $3$, the number of nodes per function has been chosen individually for each parameter between $2$ and $6$. The test has been conducted by running $10$-fold cross-validation (i.e. all predictions are done for the unseen data). For $10$ consecutive executions, the average number of the correct predictions has been observed to be $8103.7 \pm 71.0$, which gives $81.0\%\pm0.7\%$ accuracy of the model. For this example again, the accuaracy of the proposed approach is comparable to that of the Neural Designer software. The source code and the data are available online\footnote{http://ezcodesample.com/naf/reallife6.html}.

\section{Conclusions}

In this article, a novel identification method for a tree of the discrete Urysohn operators has been proposed. A particular case of such tree is the Kolmogorov-Arnold representation. The proposed algorithm can be classified as the deep machine learning algorithm, as it deals with the model consisting of several layers and having unobserved intermediate (hidden) variables. 

The suggested method is based on the fundamental properties of the discrete Urysohn operator discovered only recently. Therefore, it is a start of a new technology, which may be improved and upgraded by other researchers. The other methods, to which this technique has been compared, have been evolving for a long time, and companies, which developed the codes, have a long history in this field. Nevertheless, the proposed identification method resulted in models of a comparable accuracy on a large variety of real-life data --- a physical object (airfoil blade vibration in a wind tunnel), a biological classification problem (edible/poisonous mushrooms), an electronic dynamic system (a circuit board processing signals in a non-linear way), and a social system (human decisions regarding a subscription to a bank service).

The number of the tests, conducted by the authors, using publicly available datasets is higher than provided in this article. All tests showed same or better results compared to results reported by other researchers, when different methods have been used.

The authors can point to two major advantages of the Urysohn-tree model. First, it covers a wide range of complexity, starting from the linear regression and including the Kolmogorov-Arnold representation. Furthermore, it allows making a choice in terms of a trade-off between the accuracy and the number of parameters, preventing the overfitting. Second, it may provide a more intuitive model compared to, for example, random forest, neural networks, support vector machine, etc. Each identified function of each Urysohn operator can be shown graphically and, in some cases, the researchers can even make theoretical conclusions regarding an influence of a certain input on the output and can possibly perform a manual tuning. By looking at the functions in the piecewise-linear form, the researchers can make decisions regarding the number of the linear blocks or the number of the branch operators.

Some theoretical aspects are not yet researched, but the algorithm can already be used as is. The name suggested for it is the ``Urysohn-tree identification'' rather than the ``construction of the Kolmogorov-Arnold representation'', since the representation is a particular tree with certain hierarchy and a number of inner and outer functions, while an Urysohn tree can contain any number of Urysohn objects with at least one sequential connection and at least one unobserved variable, which is the input of one operator and the output of another operator at the same time.

\begin{flushleft}
\bibliographystyle{unsrt}
\bibliography{refs}
\end{flushleft}

\end{document}